\documentstyle[12pt]{article}

\newcommand{\Proof}{{\em Proof} \hspace{0.2in}}

\newcommand{\tag}[1]{\eqno(#1)}
\newenvironment{remark}{{\bf Remark} \hspace{0.2in}}{\\}
\newtheorem{theorem}{Theorem}[section]
\newtheorem{lemma}[theorem]{Lemma}
\newtheorem{proposition}[theorem]{Proposition}
\newtheorem{corollary}[theorem]{Corollary}

\title{On the quantization of Poisson brackets}
\author{Joseph Donin \\
Department of Mathematics and Computer Science\\
Bar-Ilan University\\
52-900 Ramat-Gan, Israel
}

\date{}
\begin{document}
\maketitle

\def\W{\wedge}
\def\t{\bigotimes}

\newcommand{\ff}{\varphi}
\def\d{\partial}
\def\k{\mbox{\bf k}}
\def\h{\hbar}
\def\ad{\mbox{ad}\,}
\def\wh{\widehat{W}(E)}
\def\wt{W(E)}
\def\ah{\widehat{A}}
\def\at{\widetilde{A}}
\def\fb{\bar{f}}

\def\ns{\mbox{nonsingular }}
\def\na{\nabla}
\def\df{\mbox{Der}_f(A)}

\def\Ah{A_{\h}}
\def\hb{\h}
\def\Ch{\k[[\hb]]}

\begin{abstract}
In this paper we introduce two classes of Poisson brackets on algebras
(or on sheaves of algebras). We call them locally free and nonsingular
Poisson brackets. Using the Fedosov's method we prove that any
locally free nonsingular
Poisson bracket can be quantized. In particular, it follows from this that
all Poisson brackets on an arbitrary field of characteristic zero
can be quantized. The well known theorem about the quantization of
nondegenerate Poisson brackets on smooth manifolds follows from
the main result of this paper as well.
\end{abstract}

\section{Introduction}
\label{Int}
{}From an algebraic point of view to quantize a commutative associative
algebra $A$ over a field $\k$ is to deform it as an associative
algebra in such a way that the deformed algebra is noncommutative.
This concept of deformation quantization was proposed in  \cite{Li}.
In many interesting cases an algebra $A$ has no
nontrivial commutative deformations, for example,  algebras
of functions on smooth or analytic manifolds. However these algebras
admit many nonequivalent noncommutative deformations.

A general theory of deformations of associative algebras has been developed
in the fundamental works of Gerstenhaber \cite{G1}, \cite{G2},
\cite{G3}, \cite{G4} (see also  \cite{GS}).
 The  Hochschild $2$-cocycles play the role
of infinitesimal objects of such deformations, hence, the tangent space
(or deformations of order one) to a ``versal family'' of deformations
of an associative algebra $A$
is the second Hochschild cohomology space $H^2(A,A)$. In the case of
commutative algebra, it is natural to begin a ``pure'' noncommutative
deformation with a skew-symmetric Hochschild cocycle, which
gives a deformation of order one. Then,
a skew-symmetric bilinear form $A\times A\to A$ is a Hochschild
cocycle if and only if it defines a biderivation
with respect to the original multiplication, i.e. satisfies the Leibniz rule.
 This form
must also satisfy the Jacobi identity if there exists an extension
of the deformation up to order two (see Section \ref{S1}),
so, the natural initial object for  deformation (or quantization)
of a commutative algebra
$A$ is a Poisson bracket $f$ on this algebra.
Given a Poisson bracket $f$ on $A$, a natural question arises:
whether there exists a deformation of $A$ over $\k[[\h]]$
with $f$ as tangent vector?
The author does not know any example of nonquantizable Poisson
bracket on ``good'' algebras.
However, there are lots of examples of Poisson brackets admitting
quantization. First of all,
any Poisson bracket on a two-dimensional smooth manifold can be
quantized. Any nondegenerate Poisson bracket
on a smooth manifold of arbitrary dimension
can be quantized as well (see \cite{deW}, \cite{Fed}, \cite{Fed1}, \cite{OM}).
Universal enveloping algebras are examples of quantizations of
linear Poisson brackets (degenerate, of course) on polynomial algebras.
It is shown in \cite{DL} that any quadratic Poisson bracket on the polynomial
algebra of three variables can be quantized. Note that quadratic
Poisson brackets correspond to deformations of polynomial algebras
as quadratic algebras.
Quantum groups are quantizations of so-called
$R$-matrix Poisson brackets on Lie groups (see \cite{Dr}).
On the other hand, the $R$-matrix Poisson brackets on
Lie groups induce, in several ways, Poisson brackets on some
homogeneous spaces of these groups. In certain cases, these Poisson
brackets can be quantized (see \cite{DGM}, \cite{DG1}, \cite{DG3}).
Other examples of quantizations can be found in \cite{GS}.

In this paper we introduce two classes of Poisson brackets on
algebras (or on sheaves of algebras) -- the class of locally free
and its subclass of locally free nonsingular Poisson brackets
(see Section \ref{S4}).
We prove that any locally free nonsingular Poisson bracket can be
quantized. It should be noted that the class of locally free
nonsingular Poisson brackets includes nondegenerate brackets
on algebras of functions on manifolds and all Poisson brackets on fields.
The class of locally free Poisson brackets includes almost
everywhere nondegenerate Poisson brackets on a manifold. In particular,
the last class contains $R$-matrix brackets on symmetric
spaces (see \cite{DG3}. The quantization of these brackets is given
in \cite{DS}.

I would like to underline that methods of the work \cite{Fed}
have influenced this paper considerably.

\section{Poisson brackets and deformations of commutative algebras}
\label{S1}
Let $A$ be an associative algebra with unit over a field $\k$ of
characteristic zero.
We will consider deformations of $A$ over the algebra of formal
power series $\Ch$ in a variable $\hb$.

By a deformation of $A$ we mean an algebra $\Ah$ over $\Ch$
that is isomorphic to $A[[\hb]]=A\hat{\t}_{\k}\Ch$ as a $\Ch$-module and
$\Ah/\h\Ah=A$ (the symbol $\hat{\t}$ denotes the tensor product
completed in the $\hb$-adic topology). We will also denote $A$ as $A_{0}$.

If $A'_{\hb}$ is another deformation of $A$, we call the deformations
$A_{\hb}$ and $A'_{\hb}$ equivalent if there exists a $\Ch$-algebra
isomorphism $A_{\hb}\to A'_{\hb}$ which induces the identity
automorphism of $A_0$.

In other words, $\Ah$ consists of elements of the form
$$x=\sum_{i=0}^{\infty}x_i\hb^i, \quad x_i\in A.$$
The multiplication in $\Ah$ is given by a $\k$-bilinear map
$F_{\hb}:A\times A\to A[[\hb]]$ written as
$$F_{\hb}(x,y)=\sum_{i\ge 0}\hb^iF_i(x,y), \quad x,y \in A,$$
where $F_0(x,y)=xy$ is the multiplication in $A$. The terms
$F_{i},\ i>0$, are
$\k$-bilinear forms $A\times A\to A$. Associativity means
that $F=F_{\hb}$ satisfies the following equation
$$F(F(x,y),z)=F(x,F(x,y)).$$
Collecting the terms by the powers of $\hb$, we obtain
$$\sum_{i+j=n}(F_i(F_j(x,y),z)-F_i(x,F_j(y,z)))=0, \quad n\geq 0. \tag{1}$$

If the multiplication in $A'_{\hb}$ is given by a bilinear map
$F'_{\hb}:A\times A\to A[[\hb]]$, then the equivalence of $\Ah$
and $\Ah '$ can be given as a power series
$$Q=Id+\sum_{i\ge 1}\hb^i Q_i,$$
where $Q_i$ are $\k$-linear maps $Q_i : A\to A$, such that
$$F'(x,y)=Q^{-1}(F(Q(x),Q(y)). \tag{2}$$

Let us consider the element $F_1$. From equation (1) for $n=1$ we
get the following relation
$$xF_1(y,z)-F_1(xy,z)+F_1(x,yz)-F_1(x,y)z=0.  \tag{3}$$
This means that $F_1$ is a Hochschild 2-cocycle.
{}From (1) for $n=2$ we have
$$F_1(F_1(x,y),z)-F_1(x,F_1(y,z))=xF_2(y,z)-F_2(xy,z)+F_2(x,yz)-F_2(x,y)z,$$
i.e. the element $F_1(F_1(x,y),z)-F_1(x,F_1(yz))$ is 3-coboundary.
If the series
$F'$ gives an equivalent deformation of $A$, then we have from (2)
$$F'_1(x,y)-F_1(x,y)=xQ_1(y)-Q_1(xy)+Q_1(x)y,$$
i.e. $F'_1(x,y)$ and $F_1(x,y)$ are cohomologous.

Throughout the remainder of this paper we will be interested in
deformations of commutative algebras. If such is the  case, we have
\begin{proposition}
\label{Pr0.1}
Let $A$ be a commutative algebra, $F_1(x,y),\ x,y\in A$ a Hochschild
cocycle,
then
$F'_1(x,y)=F_1(y,x)$ is a Hochschild cocycle.
\end{proposition}
\Proof Straightforward computation.
\smallskip \\

Thus, in the case of  commutative algebra any Hochschild
cocycle $F_1(x,y)$ can be decomposed into two cocycles:
$F_1(x,y)=\alpha(x,y)+\beta(x,y)$, where
$\alpha(x,y)=(F_1(x,y)+F_1(y,x))/2$ and
$\beta(x,y)=(F_1(x,y)-F_1(y,x))/2$.
The cocycle $\alpha(x,y)$ is commutative, i.e.
$\alpha(x,y)=\alpha(y,x)$, therefore it determines an
infinitesimal commutative deformation of $A$. We will only consider
noncommutative deformations of a commutative algebra (or quantizations),
so we will put $\alpha(x,y)=0$. Moreover, in many cases algebra $A$ will have
no commutative deformations, and any commutative Hochschild cocycle
will be a coboundary, so in these cases any deformation is equivalent to a
deformation with skew-symmetric $F_1$, i.e. we can suppose that
$F_1(x,y)=\beta(x,y)$. Put $(x,y)=\beta(x,y)$. We have
$(x,y)=-(y,x)$.
\begin{proposition}
\label{Pr0.2}

a) A skew-symmetric bilinear map $A\t A\to A$,
$x\otimes y \mapsto (x,y)$ is a Hochschild cocycle if and only if
it obeys the Leibniz rule
$$(xy,z)=x(y,z)+y(x,z).$$

b) If $((x,y),z)-(x,(y,z))$ is a Hochschild 3-coboundary, then the Jacobi
identity holds:
$$((x,y),z)+((y,z),x)+((z,x),y)=0.$$
\end{proposition}
\Proof a) If $(x,y)$ satisfies the Leibniz rule a
straightforward computation shows that it is a Hochschild cocycle.
Let $(x,y)$ be a cocycle. Then
\begin{eqnarray*}
x(y,z)-(xy,z)+(x,yz)-(x,y)z=0,\\
y(x,z)-(yx,z)+(y,xz)-(y,x)z=0,\\
(-1)[x(z,y)-(xz,y)+(x,zy)-(x,z)y]=0.
\end{eqnarray*}
Adding these equations, we will get
the Leibniz rule. b) Suppose, there exists a 2-cochain $\{x,y\}$ such
that
$$((x,y),z)-(x,(y,z))=x\{y,z\}-\{xy,z\}+\{x,yz\}-\{x,y\}z.$$
Let us also write five similar equations for permutations of  $x,y,z$.
Let $(abc)$ be a permutation of $x,y,z$.
Multiplying the equation corresponding to $(abc)$ by
$\mbox{sign}(abc)$ and adding all these equations, we
will get the Jacobi identity. The proposition is proved.
\smallskip\\

The latter proposition shows that Poisson brackets are natural initial
infinitesimal objects for noncommutative deformations of commutative
algebras. The Leibniz rule and the Jacobi identity are needed for
existence of deformations of degrees one and two in $\h$, respectively.

\section{Koszul complex and Weyl algebra}
\label{S2}
Let $A$ be a commutative associative algebra with unit over
    a field $\k$ of characteristic zero.

    Let $B$ be a bigraded
    $A$-algebra (noncommutative). We will regard $B$ as a super-algebra,
    an element of which $x\in B$ of degree $(p,q)$
    is even (odd) if the number $q$ is even (odd). Denote by $\tilde{x}$
    the parity of $x$. Then the commutator of two elements $x,y\in B$
    is defined as $[x,y]=xy-(-1)^{\tilde{x}\tilde{y}}yx$.

    An $A$-linear operator in $B$ is said to be of degree $(r,s)$
    if it sends elements of degree $(p,q)$ into elements of degree
    $(p+r,q+s)$. So the set of $A$-linear operators in $B$ may also be
    considered as bigraded super-algebra: an operator of degree
    $(r,s)$ is even or odd depending on the parity of $s$.

    An operator $D$ is called a derivation of degree $\tilde{D}$
    if the following equality
    holds
    $$D(xy)=D(x)y+(-1)^{\tilde{x}\tilde{D}}xD(y).$$
    Note that all derivations form a Lie super-algebra
    with respect to (super)-commutator. In particular, if $D$ is an odd
    derivation, then $D^2$ is an even one, and the Bianchi identity holds:
    $$[D,D^2]=D^3-D^3=0.$$

    Given an $A$-module $E$, we denote by $T(E)$, $S(E)$, and $\W E$
    the tensor, symmetric, and exterior algebra over $A$, respectively.

    Suppose $u:E \to F$ is a morphism of $A$-modules. Let us define
    an operator $d_u=d$ on the bigraded $A$-module $T(E)\t \W F$
    in the following way.
    If $x=x_1\otimes\cdots\otimes x_m\otimes y_1\wedge\cdots\wedge y_n
    \in T^m(E)\t \W^n F$, we put
    $$dx=\sum_i
    x_1\otimes\cdots\otimes\hat{x_i}\otimes\cdots\otimes x_m\otimes
    u(x_i)\wedge y_1\wedge\cdots\wedge y_n.$$

    Similarly, we define an operator $\d_u=\d$ on the bigraded
    $A$-module $S(F)\t T(E)$ assigning  to
    $x=y_1\odot\cdots\odot y_m\otimes x_1\otimes\cdots\otimes x_n
    \in S^m F\t T^n E$ the element
    $$\d x=\sum_i (-1)^{i-1}y_1\odot\cdots\odot y_m\odot u(x_i)\otimes
    x_1\otimes\cdots\otimes\hat{x_i}\otimes\cdots\otimes x_n.$$

    \begin{proposition}
    \label{Pr1.1}
    The operators $d$ and $\d$ are $A$-linear derivations of degrees
    $(-1,1)$ and $(1,-1)$, respectively, on corresponding algebras
    considered as super-algebras with respect to the second degree.
    Moreover, $d^2=0,\ \d^2=0$.
    \end{proposition}
    \Proof Straightforward computation. \smallskip \\

    Given a bilinear skew-symmetric form $\ff:\W^2E\to A$,
    denote by $I$ the ideal in $T(E)$ generated by relations
    $$x\otimes y-y\otimes x-\ff(x,y)=0. \tag{1}$$
    We will call $W(E)=T(E)/I$  the Weyl algebra associated to $\ff$.
    The operator $d$ induces a derivation on $W(E)\t\W F$.
    Indeed, $d$ applied to the left side of (1) gives zero.

    Analogously, given a bilinear symmetric form $\psi:S^2(E)\to A$,
    denote by $J$ the ideal in $T(E)$ generated by relations
    $$x\otimes y+y\otimes x-\psi(x,y)=0. \tag{2}$$
    Since $\d$ applied to (2) gives zero, $\d$ induces a derivation
    on $S(F)\t C(E)$, where $C(E)=T(E)/J$ (the Clifford algebra associated
    to $\psi$).

    In particular, if $\ff=0$ and $\psi=0$ we get Koszul complexes
    $\mbox{Kos}^{\bullet}(u)=(S(E)\t\W F,d)$ and
$\mbox{Kos}_{\bullet}(u)=(S(F)\t \W E,\d)$
    (see e.g. \cite{Ill}, pp. 107-113).

    Let us introduce an independent variable $\h$ and consider the
    modules $E[\h]=E\t_{\k}\k[\h]$ and $F[\h]=F\t_{\k}\k[\h]$ over
    the algebra $A[\h]=A\t_{\k}\k[\h]$. In this case we consider
    the following relations
    $$x\otimes y-y\otimes x-\h\ff(x,y)=0, \tag{3}$$
    $$x\otimes y+y\otimes x-\h\psi(x,y)=0 \tag{4}$$
    rather than  (1) and (2).

    Let us form as above the algebras $W(E[\h])$ and $C(E[\h])$.
    If we now regard $\h$ as being of degree two, then the relations (3)
    and (4) are homogeneous and $W(E[\h])$ and $C(E[\h])$ become graded
    algebras. Moreover, $d$ and $\d$ become derivations of
    degrees $(-1,1)$ and $(1,-1)$.

    Throughout the remainder of this paper we will often drop
    the reference to $\h$ in our notations.
    Thus, we will write $E$, $A$ rather than $E[\h]$, $A[\h]$, and so on.

    Consider in more detail the case in which $u:E\to F$ is an isomorphism
    of free modules of finite rank. If such is the  case, the Koszul complexes
    $\mbox{Kos}^{\bullet}(u)$ and $\mbox{Kos}_{\bullet}(u^{-1})$ are equal to
    $S(E)\t\W F$ as algebras
    and we get the derivations $d$ of degree $(-1,1)$ and $\d$ of degree
    $(1,-1)$ on this algebra.
    \begin{proposition}
    \label{Pr1.2}
    Let $u:E\to F$ be an isomorphism of free modules of finite rank.
    Then $(\d d+d\d)x=(p+q)x$ for $x\in S^p(E)\t\W^q F$.
    \end{proposition}
    \Proof Straightforward computation. \smallskip \\
    In particular, $(\d d+ d\d)x=0$ if and only if $x\in S^0(E)\t\W^0 F=A$.

    Notice, that there is a natural embedding of $A$-modules
    $$\sigma:S(E)\t\W F\to T(E)\t \W F,$$
    generated by the embeddings $S^n(E)\to T^n(E)$ for any $n$.

    Let $\ff$ be a skew-symmetric form on $E$, and $W(E)$ the corresponding
    factor algebra of $T(E)$, as above. Denote by $\pi$ the natural
    projection $T(E)\W F\to W(E)\W F$. By the Poincar\'e-Birkhoff-Witt theorem
    the composition
    $\pi\sigma$ gives an isomorphism of $A$-modules $S(E)\t\W F\to W(E)\t\W F$.
   Due to this isomorphism, the operator $\d$ can be carried onto $W(E)\t\W F$.
    Of course, it will not be a derivation, but the relation
    $(\d d+d\d)x=(p+q)x$ for $x\in W(E)\t\W F$
    remains true.

    Given $a\in W(E)\t\W F$, we put $s(a)=(\pi\sigma)^{-1}(a)\in S(E)\t\W F$
    and call $s(a)$ the symbol of $a$. We say that $a$ has $s$-degree $(n,m)$
    if $s(a)\in (S^n(E)\t\W^m F)[\h]$.
We will also say that $a$ has $s$-degree $n$ if $s(a)\in(S^n(E)\t\W F)[\h]$.
When $a$ has $s$-degree $n$, then
    $da$ and $\d a$ have $s$-degrees $(n-1)$ and $(n+1)$, respectively.
    It follows from the fact that the operators $d$ and $\d$ commute
    with $\pi\sigma$ and from the explicit forms of $d$ and $\d$ in
$S(E)\t\W F$.

    So, we get
    \begin{proposition}
    \label{Pr1.3}
    Let $u:E\to F$ be an isomorphism of free $A$-modules of finite rank.
    Given a skew-symmetric form $\ff$ on $E$, suppose $W(E)=W(E[\h])$ is
    the corresponding graded algebra. Then there exist a derivation $d$
    on $W(E)\t\W F$ of degree $(-1,1)$ and an $A$-linear operator $\d$ of
    degree $(1,-1)$ such that the following equality holds
    $$ (\d d+d\d)x=(p+q)x\ \mbox{  if  } s(x)\in W^p(E)\t\W^qF.  \tag{5}$$
    Moreover, $d$ and $\d$ have $s$-degree $(-1,1)$ and $(1,-1)$, respectively.
    \end{proposition}

    Let us consider the complex
$$W^{\bullet}\ :W(E)\stackrel{d}{\longrightarrow} W(E)\t\W^1 F\stackrel{d}
    {\longrightarrow}
    W(E)\t\W^2 F\stackrel{d}{\longrightarrow}\cdots. $$
    \begin{corollary}
    \label{Cor1.1}
    In the hypothesis of Proposition \ref{Pr1.3} we have
    $$H^0(W^{\bullet})=A,\ \ H^i(W^{\bullet})=0 \mbox{  for } i>0.$$
    \end{corollary}
    \Proof It follows from equality (5). \smallskip \\
    \begin{remark}
    Let us assume that $A$ is a sheaf of algebras over a topological space
    $M$, $u:E\to F$ is an isomorphism of locally free sheaves of $A$-modules,
  and $\ff$ is a mapping of sheaves $\W^2E\to A$. Then all constructions above
    make sense in this situation, including Proposition \ref{Pr1.3}
    and corollary, where $W(E)\t\W F$ becomes a sheaf of $A$-algebras.
    \end{remark}

    \section{Completed Weyl algebra and derivations}
\label{S3}
    Suppose $E$ is a free $A$-module of  finite rank,
    $\ff$ is a skew-symmetric form on $E$. Let us construct the
    Weyl algebra $W(E)=W(E[\h])$ with respect to $\ff$, as in the
    preceding section. This algebra is a graded algebra, in which
    the elements from $E$ have degree one and $\h$ is of degree two.

    The gradation in $W(E)$ induces a decreasing filtration in this
   algebra with submodules $W_p (E)$, generated by elements of degree $\geq p$.
    A nonhomogeneous element $a\in W(E)$ will be called of degree $p$,
    if $p$ is the maximal number such that $a\in W_p (E)$.

    Let us form the completion $\wh$ of $W(E)$ with respect to
    this filtration. 
    In this case $\wh$ can be regarded as a module
    over the algebra $\ah=A[[\h]]$ of power series in $\h$.
    %
    %

    Let $u:E\to F$ be an isomorphism of modules.
    It is clear that one can complete the complex
$W^{\bullet}=(W(E)\t \W F,d)$
    defined in the preceding section
    and get the complex
    $\widehat{W}^{\bullet}=(\wh\t\W F,d)$.
    Moreover, the operator $\d$ makes sense in this case and Propositions
    \ref{Pr1.2} and \ref{Pr1.3} remain true.
Further we will use the notations $W(E)$ and $W^{\bullet}$ for the completed
algebras $\widehat{W}(E)$ and $\widehat{W}^{\bullet}$.

As in Section 3 we will consider $W(E)\t\W F$ as a bigraded algebra
and as a super-algebra: an element of this algebra of the degree
$(p,q)$ is even or odd depending on the parity of $q$.
Thus, by derivations of $W(E)\t\W F$ we mean super-derivations,
by commutator of two elements from $W(E)\t\W F$ we mean
super-commutator, and so on.

    \def\fc{\widetilde{W}_0^{\bullet}(E)}

    Later we will need the following description of $A$-linear
    derivations of nondegenerate Weyl algebras.
    Assume $\ff$ is a nondegenerate skew-symmertic form on $E$.
    This means that the mapping $\phi:E\to E^\ast,\ x\mapsto \ff(x,\cdot)$,
    from $E$ into the dual module is an isomorphism.
     The Weyl algebra $W(E)$ associated to such a  $\ff$ will be called
   nondegenerate.

    \begin{proposition}
    \label{Pr2.1}
    Any $A$-linear derivation $D$ of a nondegenerate Weyl algebra $W(E)$
    is an inner one, i.e. there exists an element $v\in W(E)$ such that
    $D(x)=\frac{1}{\h}[v,x]$ for any $x\in W(E)$.
    (Recall, that we use the form $\h\ff$ in the definition of $W(E)$).
    \end{proposition}
    \Proof Let us consider the complex $(\wt\t \W E^*,d)$ associated to
    the mapping $\phi:E\to E^*,\, x\mapsto \ff(\cdot, x)$.
Let $e_i$ be a basis in $E$ over $A$
    and $e^i$ be the dual basis in $E^*$, $e^i(e_j)=\delta^i_j$.
    It is easy to verify that the operator $d$ has the form
$\frac{1}{\h}[\bar{d},\cdot]$, where
    $$\bar{d}=\sum_ie_i\otimes e^i \in E \t \W^1 E^*.$$
    One has the equalities
    $$[De_i,e_j]+[e_i,De_j]=0 \mbox{   for any } i,j.$$
    Form the element
    $$\bar{D}=\frac{1}{\h}\sum_i De_i\otimes e^i \in \wt\t\W^1 E^*.$$
    The last equation implies that $d\bar{D}=0$.
    It follows from the exactness of the complex
    $(\wt\t \W E^*,d)$ that there exists an element $v\in\wt$
    such that
    $$-dv=\frac{1}{\h}\sum_i[v,e_i]\otimes e^i=\frac{1}{\h}\sum_iDe_i\otimes
e^i,$$
\def\fh{\frac{1}{\h}}
i.e. $\fh[v,e_i]=De_i$ for all $i$. Thus the operator
   $\fh[v,\cdot]$ coincides with $D$.
The proposition is proved.
   \smallskip \\

   Let us consider derivations in the superalgebra $\wt\t\W F$ which
   (super)-commute with the multiplication by elements from $\W F$.
   Thus, the operator $d$ is such a derivation. It follows from
   the last proposition that such derivations are inner.
   \smallskip \\
   \begin{remark} As in the preceding section, all the constructions of
   this section make sense in the global case as well, when $A$ is a sheaf
   of algebras over a topological space $M$, etc. In this case
   free $A$-modules are replaced with
   locally free sheaves of $A$-modules.
   Proposition \ref{Pr2.1} is true locally. In order for any  global derivation
   $D$ in that proposition to be inner, one needs to suppose that
   $H^1(A,M)=0$, where $H^i({\cal F},M)$
denotes the $i$-th cohomology
   of a sheaf ${\cal F}$ over $M$.
   \end{remark}

   \section{Poisson brackets and Hamiltonian derivations}
\label{S4}
   We recall that a Poisson bracket on a $\k$-algebra $A$ is a skew-symmetric
   $\k$-linear form $f:\W^2A\to A$ which satisfies two conditions:\\
   the Leibniz rule
   $$f(ab,c)=af(b,c)+bf(a,c),$$
   and the Jacobi identity
   $$f(f(a,b),c)+f(f(b,c),a)+f(f(c,a),b)=0.$$
   It follows from the Leibniz rule that $f$ defines a mapping
   $\fb:A\to \mbox{Der}(A)$, namely, $\fb (a) = f(a, \cdot)$
where $\mbox{Der}(A)$ denotes the $A$-module
   of derivations of the algebra $A$.

   Note that the Jacobi identity implies that $H=\mbox{Im}(\fb)$
   forms a Lie algebra over $\k$. Indeed, let $a,b,c\in A$ and
   $\bar{a}=\fb(a), \bar{b}=\fb(b)$. We have by definition
   $\bar{a}c=f(a,c), \bar{b}c=f(b,c)$, hence $\bar{a}\bar{b}c=f(a,f(b,c))$.
   Similarly, $\bar{b}\bar{a}c=f(b,f(a,c))$. Therefore,
   $$(\bar{a}\bar{b}-\bar{b}\bar{a})c=f(a,f(b,c))-f(b,f(a,c))=f(f(a,b),c)$$
   by the Jacobi identity. This means that
   $$\bar{a}\bar{b}-\bar{b}\bar{a}=\overline{f(a,b)}, $$
   i.e. $\bar{a}\bar{b}-\bar{b}\bar{a}$ is the image of $f(a,b)$ by
   the mapping $\bar{f}$.

   It should be noted that $H$ is not an $A$-module.
   We will call elements from $H$  strong Hamiltonian derivations.
   The bracket $f$ forms a $\k$-linear nondegenerate skew-symmetric
   form $\ff$ on $H$ with
   values in $A$ in the following way. Let $x=\fb(a),\,y=\fb(b)$, then
   we put $\ff(x,y)=f(a,b)$. It is clear that such a definition of $\ff$
   makes sense.
Indeed, suppose $x_i\in A$, and $\bar{x}_i$ are the corresponding strong
Hamiltonian derivations. We must show that if $D=\sum_i
a_i\bar{x}_i=0,\,a_i\in A,$ then $\sum_i a_i\ff(\bar{x}_i,\bar{b})=0$
for any $b\in A$. But that follows from the following chain of
equalities
$$\sum_i a_i\ff(\bar{x}_i,\bar{b})=\sum_i a_i f(x_i,b)=Db=0.$$

Let us denote
   by $E$ the $A$-submodule in $\mbox{Der}(A)$ generated
   by $H$ and call the elements of it  weak Hamiltonian derivations.
   It is easy to see that $E$ forms a Lie subalgebra in
   $\mbox{Der}(A)$ over $A$, and the form $\ff$ can be extended by linearity
   to $E$.
   So, we obtain on $E$ an $A$-linear skew-symmetric form
   $\ff:\W^2 E\to A$ associated to our Poisson bracket on $A$.
We assume further that $E$ is an $A$-module of finite type.

   Note that the mapping $\phi:E\to E^*,\ x\mapsto \ff(x,\cdot)$,
   is always a monomorphism.
This follows from the fact that if $\ff(x,y)=0$ for a fixed $x$ and
any $y$, $x,y\in E$, then $x=0$. It is sufficient to prove the fact
for all $y$ of the form $\bar b$, where $b\in A$. But, if
$x=\sum_i a_i\bar{x}_i,\,a_i\in A,$ is zero, where $\bar{x}_i$
are strong Hamiltonian derivations corresponding to $x_i\in A$, then
$$\ff(x,\bar{b})=\sum_i a_i f(x_i,b)=x(b),$$
for any $b\in A$; i.e. the derivation $x$ applied to $b$ is equal to zero.

The module $E$ need not be free (or locally
   free in the global situation). We will call a Poisson bracket  free
   (locally free in the global case) if the associated  module (or sheaf)
   of weak Hamiltonian derivations $E$ is (locally) free. We will call
   a Poisson bracket $f$ \ns if the corresponding map $\phi:E\to E^*$
   is an isomorphism.

   Given a locally free Poisson bracket $f$ on $A$, let us construct
   the Weyl algebra $W(E)$ and the algebra $W^{\bullet}(E)=(W(E)\t\W
E^*,d)$ associated
   to the monomorphism $\phi:E\to E^*$, as in preceding sections.

   The Lie algebra $E$ acts on $A$ as derivations, so one can apply
   to any $a\in A$ the differential form $\nabla a,\ \nabla a(x)=xa$.
   Hence, $\nabla$ can be considered as a mapping $A\to \W^1 E^*$ with
   the property $\na(ab)=a\na(b)+b\na(a)$.

Moreover, the operator $\na$ can be extended  to  a
derivation on  the exterior algebra $\W E^*$ in the following way.
Consider $\W^n E^*$ as the algebra of $A$-linear skew-symmetric
functions on $E$ of $n$ variables with values in $A$. Then
for $\rho\in \W^{n-1} E^*$ we define $\na(\rho)\in \W^{n}E^*$ by
\begin{eqnarray*}
\na(\rho)(x_1,...,x_{n})=\sum_{1\leq i<j\leq n} & (-1)^{i+j}
\rho([x_i,x_j],x_1,...,\hat{x}_i,...,\hat{x}_j,...,x_n) \\
&+{\displaystyle \sum_{1\leq i\leq n}}
(-1)^{i-1} x_i\rho(x_1,...,\hat{x}_i,...,x_n).
\end{eqnarray*}
So the definition of $\na$ is the same as for the de Rham complex
and for the cohomology complex for Lie algebra.
It is easy to verify that $\na$ is an odd derivation of the algebra $\W E^*$
and $\na^2=0$ on $\W E^*$.
It turns out that $\na$ can be extended to the whole algebra
$W^{\bullet}=W(E)\t\W E^*$.
   \begin{proposition}
   \label{Pr3.1}
   Let $f$ be a locally free Poisson bracket on $A$. Then
   the derivation $\na$ can be extended to a derivation on $W^{\bullet}(E)$
   of degree $(0,1)$ with the property
$$d\psi=0, \tag{1}$$
where $\psi\in E\t\W^2E^*$ is the tensor of torsion defined as
$\psi(x,y)=\na_x(y)-\na_y(x)-[x,y]$
   for any $x,y\in E$ (here $[\cdot,\cdot]$ denotes the Lie bracket in $E$).
   \end{proposition}
   \Proof First of all, reduce the proposition to the case in which $E$
admits a basis over $A$ consisting of strong Hamiltonian derivations.

Let $\mbox{Spec}(A)$ be the spectrum of $A$ with Zariski topology
and structure sheaf ${\cal A}$, i.e. the affine sheme associated to $A$.
The sheaf ${\cal A}$ is a sheaf of local algebras such that
$A=H^0({\cal A})$. The Poisson bracket $f$ induces in  a natural way
the Poisson bracket $\tilde{f}$ on ${\cal A}$.
On the other hand, there is the sheaf ${\cal E}$ on $\mbox{Spec}(A)$
corresponding to $E$, which is a locally free ${\cal A}$-module,
and $E=H^0({\cal E})$.
It is easy to see that ${\cal E}$ will be the sheaf of germs of weak
Hamiltonian derivations of $\tilde{f}$.
 Denote by ${\cal E}_p$ and ${\cal A}_p$ the stalks
of the sheaves ${\cal E}$ and ${\cal A}$ over the point
$p\in \mbox{Spec}(A)$.
By definition, ${\cal E}_p$ is generated by strong Hamiltonian derivations
as ${\cal A}_p$-module. Since ${\cal E}_p$ is a free module over
the local algebra ${\cal A}_p$, one can extract from these Hamiltonian
derivations a finite basis of ${\cal E}_p$ over ${\cal A}_p$. This
basis will be a basis of the $A_U$-module $E_U$,
where $A_U$ and $E_U$ denote the spaces of sections of ${\cal A}$ and
${\cal E}$ over some neighborhood $U$ of $p$. So, we see that any
point of $\mbox{Spec}(A)$ has a neighborhood in which the module
of weak Hamiltonian derivations admits a basis consisting of strong
Hamiltonian derivations. Consider a covering ${U_i}$ of $M=\mbox{Spec}(A)$
by such neighborhoods and suppose that the proposition is true for
the algebras $A_{U_i}$ with induced Poisson brackets $f_{U_i}$.
Now, let us apply the arguments from the proof of Proposition \ref{Pr3.2}
below to the Poisson bracket
$\tilde{f}$ given on the sheaf ${\cal A}$ over the space
$M=\mbox{Spec}(A)$. Note that in this case $H^1({\cal T},M)=0$,
because ${\cal T}$ is a quasi-coherent sheaf on affine algebraic space $M$
(see the definition of ${\cal T}$ just before Proposition \ref{Pr3.2}).
So we conclude that the proposition will be true if it is true for each pair
$A_{U_i}$, $f_{U_i}$.
Thus, it suffices to prove the existence of the required derivation $\nabla$
in supposing that $A$-module $E$ admits a strong Hamiltonian basis.

 Let $e_i$ be a basis in $E$ over $A$ consisting of strong
   Hamiltonian derivations. Let us define
$\na_{e_i}(e_j)=[e_i,e_j]$
   (here we mean by
   $[\cdot,\cdot]$ the bracket in the Lie algebra $E$).
   If $x=\sum_i a_i e_i,\, a_i\in A$, is an element from $E$ we set
   $\na_{e_j}(x)=\na_{e_j}(a_i)e_i+a_i\na_{e_j}(e_i)$.
   If $y$ is another
   element from $E$, it can also be presented uniquely as a linear
   combination $\sum_ic_ie_i$ with $c_i\in A$. Hence, defining
   $\na_y(x)=\sum_ic_i\na_{e_i}(x)$, we get a linear mapping
   $\na:E\to E\t \W^1 E^*,\ x\mapsto \na.(x)$.
   The operator $\na$ defined just above is a connection along weak
   Hamiltonian derivations.

   Let us prove now that the form $\ff$ on $E$ associated to the Poisson
   bracket $f$
   is an invariant one relative to connection $\na$, i.e.
   $$\na_z(\ff(x,y))=\ff(\na_z x,y)+\ff(x,\na_z y) \mbox{  for any }x,y,z\in
E.$$
   It is sufficient to verify it for the elements $e_i$ of our basis.
   Suppose $a_i\in A$ are such elements that $e_i=f(a_i,\cdot),$
   Then, if we take $x=e_i, y=e_j, z=e_k$, the previous equation
   will be equivalent
   to the relation $f(a_k,f(a_i,a_j))=f(f(a_k,a_i),a_j)+f(a_i,f(a_k,a_j))$,
   which is true by the Jacobi identity.

   Now we can extend $\na$ on the whole algebra $T(E)\t\W E^*$ as an odd
   operator by the Leibniz rule.

Let us prove that the ideal $I\subset T(E)$ generated by the relations
$$x\otimes y-y\otimes x=\h\ff(x,y),\ \ \ x,y\in E, \tag{2}$$
is invariant under the action of $\na$.
Applying $\na$ to the left side of (2) we obtain
$$\na[x,y]=[\na x,y]+[x,\na y], \tag{3}$$
where $[x,y]=x\otimes y-y\otimes x$, the commutator in $T(E)$.
Applying $\na$ to the right side of (2) we obtain
$$\h\na(\ff(x,y))= \h(\ff(\na x,y)+\ff(x,\na y)) \tag{4} $$
because of the invariance of $\ff$.
It is obvious that the difference between (3) and (4) belongs
to $I$, which proves the invariance of $I$ under the action of $\na$.
It follows from this that $\na$ induces a well-defined derivation
on $W^{\bullet}(E)=T(E)/I\t\W E^*$.

   To verify (1), note that for elements of the basis we have
   $$\psi(e_i,e_j)=\na_{e_i}(e_j)-\na_{e_j}(e_i)-[e_i,e_j]=[e_i,e_j]$$
   by definition. On the other hand, the element $\psi$ is
   obviously $A$-bilinear.
Hence, $\psi$ has the form $\psi=\sum_{i,j}[e_i,e_j]e^i\wedge e^j$.
Recall that the operator $d$ has the form ${1 \over \h}[\bar{d},\cdot]$
(see the prooof of Proposition 4.1). Using (2) we get
$d\psi(e_i,e_j,e_k)=\mbox{Alt}_{ijk}\ff(e_i,[e_j,e_k])=
f(a_i,f(a_j,a_k))+f(a_j,f(a_k,a_i))+f(a_k,f(a_i,a_j))=0$,
by the Jacobi identity.
It implies
   that this equality holds for any elements $x,y\in E$.
   This completes the proof.\smallskip \\

   Now we want to extend the last proposition to the global situation.
   Thus, let $A$ be a sheaf of algebras over a topological space $M$
   endowed with a Poisson bracket $f$, which is a global section
   of the sheaf $Hom(\W A^2,A)$. We construct as above the sheaf
   $E$ of weak Hamiltonian derivations with the $A$-linear form $\ff$.
   The $\k$-linear differential operator $\na$ is defined on $A$ as well.

   Let us denote by $sp(E)$ the sheaf of germs of symplectic $A$-linear
   endomorphisms of $E$. By definition, such an endomorphism $Q$
   preserves the form $\ff$, i.e. $\ff(Qx,y)+\ff(x,Qy)=0,\ x,y\in E$.
   Denote by ${\cal T}$ the subsheaf of sections $s$ of
   $sp(E)\t E^*$ such that $ds=0$ for any $x,y\in E$.
   (Here we consider $sp(E)\t E^*$ as a subsheaf of $E\t E^*\wedge E^*$.)
Since $d$ is $A$-limear operator, ${\cal T}$ is a sheaf of $A$-modules.
\begin{proposition}
   \label{Pr3.2}
   Given a locally free Poisson bracket on a sheaf of algebras
   over a topological space $M$, suppose that $H^1({\cal T},M)=0$.
   Then the operator $\na$ can be extended to a derivation on $W^{\bullet}(E)$
   of degree $(0,1)$ with the property (1).
   \end{proposition}
   \Proof Let $\{U_i\}$ be an open covering of $M$ such that the sheaf $E$
   is free over each $U_i$. By Proposition \ref{Pr3.1}, there exist
   extensions $\na_i$ of $\na$ over each $U_i$. A direct check shows
   that the differences $\na_{ij}=\na_i-\na_j$ form $A$-linear
   derivations of degree $(0,1)$ over $U_{ij}=U_i\bigcap U_j$,
   and are sections of the sheaf ${\cal T}$ over $U_{ij}$.
   Moreover, they obviously form  a \v{C}ech cocycle on $M$.
   The condition $H^1({\cal T},M)=0$ means that there exist
   sections $s_i$ of ${\cal T}$ over $U_{i}$ such that
   $s_i-s_j=\na_{ij}=\na_i-\na_j$. Hence, $\na_i-s_i=\na_j-s_j$ is a global
   operator $\na '$ on $A$. It is easy to see that the operator $\na '$
   is  the required derivation on $W^{\bullet}(E)$. The proposition is proved.
   \smallskip \\

   The operator $\na$ constructed above  is a connection along weak
   Hamiltonian derivations. This connection also determines  derivations
   on the algebras $T(E)\t\W E$ and $S(E)\t\W E$ in the same way.
   Moreover, this
   connection commutes with taking of symbol, i.e.
   $$s(\na(a))=\na(s(a))\ \mbox{  for any } a\in W(E)\t\W E^*.$$

   It is easy to verify that the operators $\na^2$ and $\na d+d\na$
   are $A$-linear derivations on $W(E)\t\W E^*$
(super)commuting with elements from $\W E^*$. Let us assume that
   $\ff$ is a nondegenerate form. Then, as was shown in the
   preceding section, there exist such elements $\alpha,\beta\in W(E)
   \t\W^2 E^*$ that $\na^2$ and $\na d+d\na$ have the forms
   $\frac{1}{\h}\ad\alpha=\frac{1}{\h}[\alpha,\cdot]$
   and
   $\frac{1}{\h}\ad\beta=\frac{1}{\h}[\beta,\cdot]$,
   respectively.

   We will need the following
   \begin{lemma}
   \label{Le3.1}
   Let $\na$ satisfy the property (1). Then the elements $\alpha$ and $\beta$
   satisfy the following relations:

   $\na\alpha=0$,

   $d\beta=0$,

   $(d+\na)(\alpha+\beta)=0$.
   \end{lemma}

   \Proof There are no difficulties in seeing that the operator $\na$ preserves
   $s$-degree and that
   the elements $\alpha$
   and $\beta$ have to be of $s$-degree $(2,2)$ and $(1,2)$. By the
Jacobi identity,
   $\ad(\na\alpha)=[\na,\na^2]=0$. It follows from this that $\na\alpha$
   must commute with all elements from $W(E)$, i.e. to be of $s$-degree zero.
   But on the other hand, $\na\alpha$ has $s$-degree $2$. This implies that
   $\na\alpha=0$, because the center of the algebra $W(E)$ consists of
elements of $s$-degree zero, which proves the first relation.

   Further, the operator $d$ is realized as $\frac{1}{\h}[\bar{d},\cdot]$
   where $\bar{d}=\sum_i e_i\otimes e^i$, so
   $\ad(\beta)=\ad(\na\bar{d})$.
   This means that $\beta=\na\bar{d}+c$, where $c$ is some central element.
   It implies that $d\beta=d(\na\bar{d})$.
But it is easy to check that
   $\na\bar{d}=(\na_{e_i}e_j-\na_{e_j}e_i-[e_i,e_j])\otimes e^i\W e^j=\psi$.
So, the second relation follows from the property (1).

   To prove the third relation, note that $a=(d+\na)(\alpha+\beta)=
   d\alpha+\na\beta$ due to the first two relations. On the other hand,
   $\ad a=[d+\na,(d+\na)^2]=0$ by the Bianchi identity (see \ref{S1}).
   But the $s$-degree
   of $a$ has to be equal to 1, because two summands of $a$ have this
   $s$-degree. So, $a$ must be equal to zero, which proves the third
   relation of the lemma.

   \section{A topological lemma}
\label{S5}
   In this section we will denote by
    $E$ an Abelian group with filtration $E=\cdots\supset E_i\cdots$.
   We will assume that $i$ runs over all the integers and $\bigcup E_i=E$,
   $\bigcap E_i=0$. The degree of an element  $x\in E$ is the maximal
   number $i$ such that $x\in E_i$. We will suppose that any element from
   $E$ has a finite degree, and denote by $deg(x)$ the degree of $x$.

   The filtration defines on $E$ a topology: the  $E_i$ form a fundamental
   system of neighborhoods of zero. Every  group with such a filtration
  can be completed with respect to that topology. Henceforth we assume that $E$
   is complete.

   Let $\Phi: E\to E$ be an arbitrary mapping (in the set-theoretic sense).
   The following simple lemma gives a criteria for the operator $Id+\Phi$
   to be invertible.
   \begin{lemma}
   \label{Le4.1}
   Let $E$ be a complete Abelian group with filtration.
   Suppose that an operator $\Phi: E\to E$ satisfies the following condition:
   $$deg(\Phi(x)-\Phi(y))>deg(x-y).$$
   Then the operator $Id+\Phi$ is invertible.
   \end{lemma}
   \Proof The lemma will be proved if we establish the existence and uniqueness
   of a solution of the equation $b=x+\Phi(x)$, where $b\in E$ is given.
   Put
   $$x_0=b-\Phi(b),\ x_1=b-\Phi(x_0),\ ...\ ,x_k=b-\Phi(x_{k-1}),\ ...$$
   The sequence $(x_k)$ is convergent. Indeed,
   $x_{k+1}-x_k=\Phi(x_k)-\Phi(x_{k-1})$. This implies that
   $$deg(x_{k+1}-x_k)=deg(\Phi(x_k)-\Phi(x_{k-1}))>deg(x_k-x_{k-1}),$$
   which proves the convergence.
   Let $x$ be the  limit of this sequence. Then it is obvious that $x$
   is a solution of our equation. If $x'$ is an other solution, then
   it should be $deg(x-x')=deg(\Phi(x')-\Phi(x))>deg(x'-x)$,
   which is impossible. This proves uniqueness and the Lemma.\smallskip \\
   There are no difficulties in seeing that both operators $Id+\Phi$
and its inverse are continuous in the topology associated to the filtration.

   \section{Complex associated to Poisson bracket}
\label{S6}
   In this section we will again suppose that $A$ is an algebra with
   a Poisson bracket $f$, and $E$ is the module of weak Hamiltonian
   derivations of $A$. We also suppose that $f$ is locally free and
nonsingular.
   This means that the induced morphism $\phi:E\to E^*$ is an isomorphism
   (see \ref{S3}). So, we can construct the completed complex
    $$W^{\bullet}(E)=(W(E)\t\W E^*,d),$$
   as in \ref{S3}.
   Recall that the operator $\d$ is well defined on this complex.
   It is clear that the operator $\na$ constructed in \ref{S4}
   being of degree $(0,1)$ is well defined in this complex as well.
Denote by $W_i(E)$ the set of elements from $W(E)$ of degree $\geq 2$.

\def\Wd{W^{\bullet}(E)}
   Given $r\in W^{\bullet}(E)$, let us denote by $\ad r$ the
   inner derivation $[r,\cdot]$ in $\Wd$, where
   the bracket $[\cdot,\cdot]$ is regarded as a commutator in superalgebra.

   We want to construct such an element $r\in W_2(E)\t\W^1 E^*$
   that the derivation $D=d+\na+\fh\ad r$ will satisfy the property
   $D^2=0$. Let us find such $r$.

   We have
   $$D^2=\na^2+d\na+\na d+\fh[d,\ad r]+\fh[\na,\ad r]+(\fh\ad r)^2. \tag{1}$$
   Using the fact that $\na,d,\ad r$ are odd, and using the Jacobi
identity in the super-case, we get
that $[d,\ad r]=\ad(dr),[\na,\ad r]=\ad(\na r),
   (\ad r)^2=\frac{1}{2}\ad [r,r]$.
   Moreover, a direct computation shows that the derivation
   $B=\na^2+d\na+\na d$ is $A$-linear. So, by Proposition \ref{Pr2.1}
   (see also Lemma \ref{Le3.1}) there exists an element
   $b=\alpha+\beta\in W_1(E)\t\W^2 E^*$
    such that $B=\frac{1}{\h}[b,\cdot]=\frac{1}{\h}\ad b$.
   So we can rewrite (1) in the form
   $$D^2=\ad(\frac{1}{\h}b+\fh dr+\fh\na r+\frac{1}{2}[\fh r,\fh r]). \tag{2}$$
   Therefore, if we will be able to find an element $r$ such
   that the equation
   $$\frac{1}{\h}b+\fh dr+\fh\na r+\frac{1}{2}[\fh r,\fh r]=0 \tag{3}$$
   holds, then the condition $D^2=0$ will be satisfied  with this $r$.
   \begin{proposition}
   \label{Pr5.1}
   There exists  $r\in W_2(E)\t\W^1 E^*$ such that
   equation (3) holds.
   \end{proposition}
   \Proof At first, let us slightly modify the operator $\d$.
   We set $\delta(a)=\frac{1}{p+q}\d(a)$ for  elements
   $a\in \Wd$ of $s$-degree $(p,q)$. Therefore, when $a$ is of $s$-degree
$(p,q)$
   with $q>0$ we have
   $$(d\delta+\delta d)a=a. \tag{4}$$
   Of course, $\delta^2=0$ as well.

   The element $-\delta(b)$ obviously belongs to
   $W_2(E)\t\W^1 E^*$. Let us consider the following
   equation in $W_2(E)\t\W^1E^*$
   $$-\delta(b)=r+\delta(\na r +\frac{1}{2\h}[r,r]). \tag{5}$$
   It is easy to see that the operator
$\Phi(r)=\delta(\na r +\frac{1}{2\h}[r,r])$
   satisfies the hypothesis of Lemma \ref{Le4.1}, because $\delta$
   increases degree by one  with respect to the filtration defined on $\wt$.
   So, using this lemma we can find  $r\in W_2(E)\t\W^1E^*$
   satisfying (5).

   Let us show that this $r$ satisfies (3) as well.
   Denote by $a$ the left side of (3).
   Note that $\delta r=0$ by (5), so $\delta dr=r$
   by (4). It implies that $\delta a=0$.

   A direct computation shows that
   $$(d+\na+\fh\ad r)a=(d+\na)(\frac{1}{\h}b)=\fh(d+\na)(\alpha+\beta).$$
   But the right-hand side expression is equal to zero due to
Lemma \ref{Le3.1}.
   So,
   $$Da=da+\na a+\fh[r,a]=0. \tag{6}$$
   Since $\delta a=0,\ \delta da=a$ by (4), and we get from (6)
   $a=-\delta(\na a +\fh[r,a])$. But this is possible only if $a=0$, because
   the degree of the right side of this equation is greater at least by one
   than the degree of the left side. It proves that $r$ satisfies the
equation (3).
   The proposition is proved.\smallskip \\

   So, we have constructed a derivation on the algebra $\Wd$
   of the form $D=d+\na+\fh\ad r$ and such that $D^2=0$.
   Therefore, one can write the following complex
   $$DW^{\bullet} :W(E)\stackrel{D}{\longrightarrow}
   W(E)\t\W^1 E^*\stackrel{D}{\longrightarrow}
   W(E)\t\W^2 E^*\stackrel{D}{\longrightarrow}\cdots. $$

   The cohomology $H^0(DW^{\bullet})$ of this complex is a subalgebra $A_{\h}$
   of the algebra $W(E)$. In the next section we will show
   that the complexes $DW^{\bullet}=(W^{\bullet},D)$
   and $(W^{\bullet},d)$ are isomorphic as $\k[[\h]]$-modules, so
   $A_{\h}$ coincides with $A[[\h]]$ as $\k[[\h]]$-module. Moreover, we will
   see that $A_{\h}$ is a quantization of $A$ by our Poisson bracket $f$.
   \smallskip \\
   \begin{remark} As in the preceding sections, note that
   in the case when $A$ is a sheaf
   of algebras over a topological space $M$ the construction of
   the derivation $D$ can be realized globally.
   \end{remark}

   \section{Quantization}
\label{S7}
   \begin{proposition}
   \label{Pr6.1}
   There exists a $\k[[\h]]$-linear operator $Q$ on $W(E)$
   such that $d=QDQ^{-1}$, therefore, this operator gives an isomorphism
   of the complexes
   $DW^{\bullet}=(\Wd,D)$
   and $(\Wd,d)$.
(We assume here that $Q$ acts on $W(E)\t\W E^*$ as $Q\t 1$.)
   \end{proposition}
   \Proof
   Let us put $Q=Id+\delta(\na+\fh\ad r)$ and prove that it is an operator
  as  required in the proposition. First of all, $Q$ is invertible by lemma
   \ref{Le4.1}, because $\delta$ increases degree.
   We have to show that $dQ-QD=0$,
   i.e.
   $$d(Id+\delta(D-d))-(Id+\delta(D-d))D=0. \tag{1}$$
   But $\delta(D-d)D=-\delta dD=-\delta d(D-d)$, because $D^2=d^2=0$.
   Using  this in (1) we get $d-D +(d\delta+\delta d)(D-d)=0$, which
   is true because $d-D$ is a derivation
of degree $1$ with respect to the second degree,
   so $(d\delta+\delta d)(D-d)=D-d$. Proposition is proved.\smallskip\\

   Thus, the subalgebras $A[[\h]]$ and $A_{\h}$ of $W(E)$
   are isomorphic as $k[[\h]]$-modules and $Q^{-1}:A[[\h]]\to A_{\h}$
   realizes this isomorphism.

   The operator $Q^{-1}$ has the form
   $$Q^{-1}=Id-\delta(\na+\fh\ad r)+(\delta(\na+\fh\ad r))^2+\cdots.$$
   Let us apply it to
   elements $a,b\in A$. We obtain
   $$Q^{-1}a=a-\delta\na(a)+\cdots,\ Q^{-1}b=b-\delta\na(b)+\cdots.$$
   Taking into account that $A$ lies
   in the center of $W(E)$, we get
   $$[Q^{-1}a,Q^{-1}b]=[\delta\na(a),\delta\na(b)]+\cdots.\tag{2}$$
   But the first bracket in the right-hand side expression is equal to
   $\h f(a,b)$, which follows from the definitions of  action of $\na$ on
   elements from $A$ and commutation in $W(E)$. Using  the
   operator $Q$, the algebra $A_{\h}$ can be identify with $A[[\h]]$ as
   $\k[[\h]]$-module and
    the new multiplication in $A$ has the form
   $a*_{\h}b=Q(Q^{-1}(a)Q^{-1}(b))$.
Taking into account that the operator $\delta(\na+\fh\ad r)$
increases $s$-degree by one  and the fact that the element $a*_{\h}b$
has $s$-degree zero, one can deduce, using (2), that the coefficient
of $\h$ in $Q(Q^{-1}a,Q^{-1}b)$ is equal to $f(a,b)$ and the other
terms have order in $\h$ greater than one.
   So, we have proved
   \begin{proposition}
   \label{Pr6.2}
   Let $A$ be a sheaf of algebras on a topological space $M$,
   $f$ a locally free \ns Poisson bracket on $A$. Suppose
   that $H^1({\cal F},M)=0$ for the sheaves
   of $A$-modules over $M$ (mentioned in preceding sections).
   Then there exists a quantization $A_{\h}$ of $A$ by the bracket $f$.
   \end{proposition}
   {\bf Remarks}

a) The construction of quantization shows that if
   $a\in A$ is an element such that $f(a,b)=0$ for any $b\in A$,
   then $a$ lies in the center of $A_{\h}$ with respect to the new
multiplication.

b) If a family $f_t$ of Poisson brackets is given, the construction
of quantization shows that this family can be quantized simultaneously.

   \begin{corollary}
   \label{Cor1}
   Let $K$ be a field of finite transcendence degree over $\k$.
   Then any Poisson bracket on $K$ can be quantized.
   \end{corollary}
   \Proof Indeed, consider $K$ as a sheaf over a point. Since $K$ is
of finite transcendence degree over
$\k$, the weak Hamiltonian derivations $E$ form a   finite-dimensional
vector space over $K$, therefore $E$ is a free $K$-module.
The mapping $\phi:E\to E^*$ induced with the Poisson bracket, being a
monomorphism (see \ref{S4}) of vector spaces of the same
dimensions, is an isomorphism. Hence, any Poisson bracket on the
field is locally free and nonsingular, by definition,
and the corollary follows from Proposition \ref{Pr6.2}.\smallskip \\
   \begin{corollary}
   \label{Cor2}
   (\cite{deW}, \cite{Mel}, \cite {Fed})
   Let $M$ be a smooth manifold. Then any nondegenerate Poisson
   bracket on $M$ can be quantized.
   \end{corollary}
   \Proof It is clear that any nondegenerate Poisson bracket on
   $M$ will be locally free and nonsingular. Moreover, $H^1({\cal F},M)=0$
   for any sheaf of modules over the sheaf of algebras of smooth
functions on $M$.
   \smallskip \\

The same argument shows that  Corollary \ref{Cor2} remains true if one
replaces the smooth manifold
   $M$ by a complex analytic Stein manifold or an affine algebraic
 smooth variety. Indeed, all the sheaves ${\cal F}$ from Proposition
\ref{Pr6.2} are coherent sheaves of modules over the structure sheaf.
 But, for these classes of spaces $H^1({\cal F},M)=0$
for any coherent sheaf of modules.

We can also formulate an assertion on the quantization
of an arbitrary Poisson bracket
on Stein analytic spaces or affine algebraic varieties.

For example, let $X$ be a reduced (i.e. without nilpotent elements in the
structure sheaf) complex analytic space,
  and $A$ the structure sheaf of $X$. Assume that ${\cal F}$ is a coherent
analytic  sheaf on $X$. It follows from the existence of a free resolution
  of ${\cal F}$ that one can find an analytic subset $Y\subset X$
  of codimension one such that ${\cal F}$ will be locally free
  over $A_Y$. Here $A_Y$ denotes the sheaf of meromorphic functions
  on $X$ with poles in $Y$ (see \cite{GR}).

  It follows from this that if $f$ is an arbitrary Poisson bracket
  on $X$, then there exists an analytic subset $Y\subset X$ of codimension
  one such that the corresponding sheaf $E$ of weak Hamiltonian
derivations is locally free, and the monomorphism $\phi:E\to E^*$ induces an
isomorphism $A_Y\t_A E \to A_Y\t_A E^*$
of $A_Y$-modules. Therefore
$f$  determines a locally free \ns Poisson bracket on $A_Y$.
  So, we obtain from Proposition \ref{Pr6.2}
  \begin{corollary}
  Let $X$ be an reduced complex analytic Stein space. Suppose $f$ is a Poisson
  bracket on $X$. Then there exists an analytic subset $Y\subset X$
  of codimension one such that the algebra $A_Y$ can be quantized by
  this bracket.
  \end{corollary}

  A similar statement is valid for affine algebraic varieties.
\smallskip \\

\noindent  {\bf Acknowledgements.} The author is happy to thank
Lenny Makar-Limanov and Steven Shnider
  for stimulating discussions and very helpful remarks.
I thank Martin Bordemann and Claudio Emmrich who noticed that my
proof of first version of Proposition 5.1 was noncorrect. This allowed me to
reformulate that Proposition.

\end{document}